\documentstyle[amssymb,12pt]{article}
\title{A nonstructure theorem for countable, stable, unsuperstable theories}
\author{Michael C.\
Laskowski\thanks{Partially supported
by NSF grant DMS-0600217}\\Department of Mathematics\\
University of Maryland
\and 
S. Shelah
\thanks{
Partially supported by U.S.-Israel Binational Science Foundation 
Grant no.\ 2002323 
and Israel Science Foundation Grant no.\ 242/03. 
 Publication no.\ 871.
}\\Department of Mathematics\\
Hebrew University of Jerusalem\\
Department of Mathematics\\
Rutgers University
}

\def\cbar{\bar{c}}

\def\xbar{\bar{x}}
\def\ybar{\bar{y}}
\def\phi{\varphi}
\def\C{{\frak  C}}

\def\Rinfty{{\bf R}^\infty}
\def\ID{{\cal  ID}}
\def\tp{{\rm tp}}
\def\stp{{\rm stp}}

\def\dom{{\rm dom}}
\def\acl{{\rm acl}}
\def\dcl{{\rm dcl}}
\def\bp{{\bf Proof.}\quad}
\def\endproof{\medskip}
\def\<{\langle}
\def\>{\rangle}
\newtheorem{Theorem}{Theorem}[section]
\newtheorem{Proposition}[Theorem]{Proposition}
\newtheorem{Definition}[Theorem]{Definition}

\newtheorem{Lemma}[Theorem]{Lemma}

\begin{document}
\maketitle

\begin{abstract}
A trichotomy theorem for countable, stable, unsuperstable
theories is offered.  We develop the notion of a `regular
ideal' of formulas and study types that are minimal
with respect to such an ideal.
\end{abstract}

\section{Introduction}
By definition, a stable unsuperstable theory admits
a type that is not based on any finite subset of its domain.
\relax From this one sees that such a theory admits trees of
definable sets.  That is, there is a sequence $\<\phi_n(x,y):n\in\omega\>$
of formulas such that for any cardinal $\kappa$ there are
definable sets $\{\phi_n(x,a_\nu):\nu\in {^{<\omega}\kappa}\}$
giving rise to $\kappa^{\aleph_0}$ partial types
$\{p_\mu:\mu\in {^\omega\kappa}\}$ where each $p_\mu$ forks over
$\{a_{\mu|k}:k<n\}$ for all $n\in\omega$.  In~\cite{Shc} 
the second author used these trees to count the number of uncountable models
or to find the maximal size of a family of pairwise nonembeddable
models of a fixed cardinality of any stable, unsuperstable theory.
However, for other combinatorial questions, such as computing
the Karp complexity of the class of uncountable models of such a theory,
the existence of these trees does not seem to be sufficient.
Here we prove that when the language is countable, any strictly
stable theory exhibits one of three more detailed nonstructural
properties.  This trichotomy is used in \cite{LSh2}, but it is likely
to be used in other contexts as well.  Two of the alternatives,
the Dimensional Order Property (DOP) or a theory being deep
appear in \cite{Shc} and are compatible with superstability.
The third alternative is new and is captured by the
following definition:

\begin{Definition}  \label{witness}
{\em  
An {\em abelian group witness to unsuperstability\/}
is a descending sequence $\<A_n:n\in\omega\>$ of abelian groups
with $[A_n:A_{n+1}]$ infinite for each $n$ such that
the intersection $A=\bigcap_n A_n$ is connected and whose
generic type is regular.
}
\end{Definition}

The existence of such a sequence readily contradicts superstability
as for any cardinal $\kappa$ one immediately obtains a tree 
$\{C_\mu:\mu\in {^\omega\kappa}\}$ of cosets of $A$.
As well, with Theorems~\ref{ss} and \ref{rperp} we see that one
can frequently say more about the generic type of $A$.
This added information is used in \cite{LSh2}.

In order to establish these results, the bulk of the paper discusses
the notion of a {\em regular ideal\/} of formulas (see Definition~\ref{r}).
The origins of these ideas date back to Section V.4 of ~\cite{Shc}
and have been reworked and expanded in \cite{Hr1} and \cite{Pillay}.

Our notation is standard, and complies with either \cite{Pillay} or
\cite{Shc}.  For a stable theory $T$  $\kappa_r(T)$
denotes the least regular cardinal $\kappa$ such that there is no
forking chain of length $\kappa$.  Thus, a stable theory is superstable
if and only if $\kappa_r(T)=\aleph_0$ and $\kappa_r(T)=\aleph_1$
when $T$ is countable and strictly stable.  We call a model
`a-saturated' (a-prime) in place of `${\bf F}^a_{\kappa_r(T)}$-saturated'
(${\bf F}^a_{\kappa_r(T)}$-prime).

{\bf Throughout the whole of this paper we assume `${\bf T=T^{{\rm eq}}}$.'}
That is, $T$ is a stable theory in a multi-sorted language,
$\C$ is a large, saturated model of $T$, and the language $L$
is closed under the following operation:
If $E(\xbar,\ybar)$ is a definable 
equivalence relation then there is a sort $U_E$ and a definable
surjection $f_E:\C^{lg(\xbar)}\rightarrow U_E(\C)$ in the language
$L$.  In particular, the set of sorts is closed under finite
products.  Thus
any finite tuple of elements from varying sorts can be viewed as an
element of the product sort.  With this identification,
every formula can be considered to have a single free variable.
As notation,  $L(\C)$ denotes the set of formulas with parameters from
$\C$ and for a specific sort $s$,
$L_s(\C)$ denotes the $L(\C)$-formulas
$\phi(x)$  in which the free variable has sort $s$.

\section{Regular ideals}

\begin{Definition} {\em  An {\em invariant ideal $\ID$\/} is a subset of 
$L(\C)$ containing all algebraic formulas, closed under
automorphisms of $\C$, and for any sort $s$ and any $\phi,\psi\in L_s(\C)$
\begin{enumerate}
\item If $\phi,\psi\in\ID$ then $\phi\vee\psi\in\ID$; and
\item If $\phi\vdash\psi$ and $\psi\in\ID$, then $\phi\in\ID$.
\end{enumerate}
A partial type $\Gamma$ (i.e., a subset of $L_s(\C)$ for some sort $s$)
is {\em $\ID$-small\/} if it entails some element of $\ID\cap L_s(\C)$.
}
\end{Definition}

Many times we will make use of the fact that formulas in $\ID$ may have
`hidden' parameters.

\begin{Lemma} \label{1}  Let $\ID$ be any invariant ideal.
\begin{enumerate}
\item  A complete type $p\in S(A)$ is $\ID$-small if and only if $p\cap\ID\neq\emptyset$.
\item  For any $A$ and $a$, $\stp(a/A)$ is $\ID$-small if and only if $\tp(a/A)$ is
$\ID$-small.
\item  If $A\subseteq B$ and $\tp(a/B)$ does not fork over $A$, then
$\tp(a/A)$ is $\ID$-small if and only if $\tp(b/A)$ is $\ID$-small.
\end{enumerate}
\end{Lemma}

\bp  (1) Right to left is immediate.  For the converse, assume 
$p$ entails $\psi\in\ID$.  By compactness there is $\phi\in p$ such
that $\phi\vdash\psi$, hence $\phi\in\ID$.
   
(2) Right to left is clear.  If $\stp(a/A)$ entails $\psi(x,b)\in\ID$,
then by compactness and the finite equivalence relation theorem
there is an $A$-definable equivalence relation $E(x,y)$
with finitely many classes
such that $\tp(a/A)\cup\{E(x,c)\}\vdash\psi(x,b)$ for some $c$. Choose 
$A$-automorphisms $\{\sigma_i:i<n\}$ of $\C$ such that
$\{E(x,\sigma_i(c)):i<n\}$ includes all the $E$-classes.
Since $\ID$ is an invariant ideal $\bigvee_{i<n}\psi(x,\sigma_i(b))\in\ID$
and $\tp(a/A)\vdash \bigvee_{i<n}\psi(x,\sigma_i(b))$.

(3)  By (2) it suffices to prove this for strong types.  Assume 
$\stp(a/B)$ is $\ID$-small.  By (1) and (2), choose
$\psi(x,b)\in \tp(a/B)\cap\ID$.  Choose $\{b_i:i\in\kappa(T)\}$ independent
over $A$, each having the same strong type over $A$ as $b$.  
Since $\ID$ is invariant, $\psi(x,b_i)\in\ID$ for each $i$.
Furthermore, 
since any $a'$ realizing $\stp(a/A)$ is independent with some $b_i$ over $A$,
$ab$ and $a'b_i$ realize the same strong type over $A$, hence
$\psi(a',b_i)$ holds.  By compactness, there is a finite subset $F$
such that $\stp(a/A)\vdash\bigvee_{i\in F} \psi(x,b_i)$, so $\stp(a/A)$ is
$\ID$-small.
\endproof

\begin{Definition}  \label{r}
{\em 
An invariant ideal $\ID$ is {\em regular\/} if, for all $L(\C)$-formulas
$\psi(y)$ and $\theta(x,y)$, {\bf IF} $\psi\in\ID$ and
$\theta(x,b)\in\ID$ for every $b\in\psi(\C)$ {\bf THEN}
$\exists y(\psi(y)\wedge\theta(x,y))\in\ID$.
}
\end{Definition}

We call a strong type $\stp(a/A)$ {\em $\ID$-internal\/} if there is
a set $B\supseteq A$ independent from $a$ over $A$, a $B$-definable
function $f$, and elements $\cbar$ such that $\tp(c/B)$ is $\ID$-small
for each $c\in\cbar$ and $a=f(\cbar)$.  The strong type $\stp(a/A)$
is {\em $\ID$-analyzable\/} if there is a finite sequence $\<a_i:i\le n\>$ from
$\dcl(Aa)$ such that $a_n=a$ and $\stp(a_i/A\cup\{a_j:j<i\})$ is
$\ID$-internal for each $i\le n$.  Since $\ID$ is a collection of formulas,
this definition of analyzability is equivalent to the usual one, see e.g.,
\cite{Pillay}.

In order to iterate the defining property of a regular ideal, we need the
following notion, whose terminology is borrowed from \cite{HrSh}.

\begin{Definition} {\em  A formula $\phi(x,c)$ 
 is {\em in $\ID$, provably over $B$\/} if there is
some $\theta(y)\in\tp(c/B)$ such that $\phi(x,c')\in\ID$ for every
$c'$ realizing $\theta$.
}
\end{Definition}

\begin{Lemma}  \label{iterate}
For all sets $B$ and every $n\in\omega$, if $\phi(x,y_0,\dots,y_{n-1})$
is $B$-definable and $a,c_0,\dots,c_{n-1}$ satisfy:
\begin{enumerate}
\item  $\tp(c_i/B)$ is $\ID$-small for each $i<n$;
\item  $\phi(x,c_0,\dots,c_{n-1})\in\ID$ provably over $B$; and
\item  $\phi(a,c_0,\dots,c_{n-1})$
\end{enumerate}
then $\tp(a/B)$ is $\ID$-small.
\end{Lemma}

\bp  Fix any set $B$.  We argue by induction on $n$.  If $n=0$ the formula
$\phi(x)$ itself witnesses that $\tp(a/B)$ is $\ID$-small.
Assume the result holds for $n$ and fix a formula $\phi(x,c_0,\dots,c_n)$
and $a,c_0,\dots,c_n$ as in the hypotheses.  Choose a formula 
$\theta(y_0,\dots,y_n)\in\tp(c_0\dots c_n/B)$ such that
$\phi(x,c_0',\dots,c_n')\in\ID$ for all $c_0'\dots c_n'$ realizing $\theta$
and, using Lemma~\ref{1}, choose $\psi(y_n)\in\tp(c_n/B)\cap\ID$.

Let $\theta^*(y_0,\dots,y_n):=\theta(y_0,\dots,y_n)\wedge\psi(y_n)$,
$\theta'(y_0,\dots,y_{n-1}):=\exists y_n\theta^*$ and
$$\phi'(x,y_0,\dots,y_{n-1}):=\exists y_n(\phi(x,y_0,\dots,y_n)\wedge\theta^*(y_0,\dots,y_n))$$
We argue that Conditions (1)--(3) are satisfied by 
$\phi'$ and $a,c_0,\dots,c_{n-1}$.
Conditions (1) and (3) are clear.  We claim that the formula
$\theta'$ witnesses that $\phi'(x,c_0,\dots,c_{n-1})\in\ID$ provably
over $B$.  Indeed, it is clear that $\theta'\in\tp(c_0\dots c_{n-1}/B)$,
so choose $c_0'\dots c_{n-1}'$ realizing $\theta'$.  
Since $\psi\in\ID$ and $\theta^*(c_0',\dots,c_{n-1}',y_n)\vdash\psi$,
$\theta^*(c_0',\dots,c_{n-1}',y_n)\in\ID$.  As well, 
for any $c_n'$ such that $\theta^*(c_0',\dots,c_n')$ holds,
we have $\theta(c_0',\dots,c_n')$ holding as well,
so $\phi(x,c_0',\dots,c_n')\in\ID$.  Thus
$\phi'(x,c_0',\dots,c_{n-1}')\in\ID$ since $\ID$ is a regular ideal.
\endproof

\begin{Proposition}  \label{internal}
If $\stp(a/A)$ is $\ID$-internal, then
$\tp(a/A)$ is $\ID$-small.
\end{Proposition}

\bp  Choose $B\supseteq A$ independent from $a$ over $A$, 
a $B$-definable formula $\phi(x,\ybar)$, and
a tuple of elements
$\cbar$ such that each $\tp(c/B)$ is $\ID$-small for each $c\in\cbar$,
$\phi(a,\cbar)$ holds, and $\exists^{= 1}x\phi(x,\cbar)$.
But the formula $\phi(x,\cbar)\in\ID$ provably over $B$ via the formula
$\exists^{=1}x\phi(x,\ybar)$, so $\tp(a/B)$ is $\ID$-small
by Lemma~\ref{iterate}.  That $\tp(a/A)$ is $\ID$-small 
follows from Lemma~\ref{1}.
\endproof

The reader is cautioned that while $\ID$-internal types are
$\ID$-small, this result does not extend to $\ID$-analyzable types.
In fact, the theory and type mentioned in Remark~8.1.6 of \cite{Pillay}
gives rise to an example of this.  Much of the motivation of
this section, and in particular how it differs from
treatments in \cite{Hr1} and \cite{Pillay}, 
revolves around  how we handle $\ID$-analyzable types that are
not $\ID$-small. 

\begin{Definition}  \label{foreign} {\em 
A strong type $p$ is {\em foreign to $\ID$,\/} written
$p\perp\ID$, if $p\perp q$ for every $\ID$-small $q$.
}
\end{Definition}

\begin{Lemma}  \label{charforeign}
The following are equivalent for any regular ideal $\ID$ and
any strong type $p$:
\begin{enumerate}
\item  $p\perp \ID$;
\item $p\perp q$ for every $\ID$-internal strong type $q$;
\item $p\perp q$ for every $\ID$-analyzable strong type $q$;
\item If $p=\stp(a/A)$ then there is no $a'\in\dcl(Aa)$
such that $\tp(a'/A)$ is $\ID$-small.
\end{enumerate}
\end{Lemma}

\bp  $(1)\Rightarrow (2)$ follows immediately from 
Proposition~\ref{internal}.  $(2)\Rightarrow(3)$
follows by induction on the length of the $\ID$-analysis,
using the fact that $p\perp \tp(b/B)$ and $p\perp \tp(a/Bb)$
implies $p\perp\tp(ab/B)$.
$(3)\Rightarrow(4)$ is trivial, and $(4)\Rightarrow(1)$
follows immediately from (say) Corollary~7.4.6 of \cite{Pillay}.
\endproof

The reader is cautioned that when the regular ideal is not
closed under $\ID$-analyzability, these definitions differ from
those in \cite{Pillay}.

\begin{Definition} {\em  A partial type $\Gamma$ is {\em $\ID$-large\/}
if it is not $\ID$-small.  $\Gamma$ is {\em $\ID$-minimal\/} if
it is $\ID$-large, but any forking extension of $\Gamma$
is $\ID$-small.  
$\Gamma$ is {\em $\ID_\perp$-minimal\/} if it is $\ID$-large, but
any forking extension $\Gamma\cup\{\theta(x,c)\}$ is $\ID$-small whenever
$\stp(c/\dom(\Gamma))\perp\ID$.
}
\end{Definition}

Clearly  $\ID$-minimality implies  $\ID_\perp$-minimality,
but one of the applications in  Section~\ref{app} will use $\ID_\perp$-minimal
types that are not $\ID$-minimal.  

\begin{Lemma}  \label{regular} Let $\ID$ be any regular ideal.  
If a strong type $p$ is both $\ID_\perp$-minimal and foreign to $\ID$,
then $p$ is regular.
\end{Lemma}

\bp  The point is that a counterexample to the regularity of
$p$ can be found within the set of realizations of $p$.
If $M$ is a-saturated and
$p=\tp(a/M)$ is not regular then there are 
a tuple
$\cbar=\<c_1,\dots,c_n\>$ realizing $p^{(n)}$
for some $n$ and a realization $b$ of $p$
such that $\tp(a/M\cbar)$ forks over $M$,
$\tp(b/M\cbar)$ does not fork over $M$,
and $\tp(b/M\cbar a)$ forks over $M\cbar$.
Let $q=\tp(a/M\cbar)$ and
choose an $L(M)$-formula $\theta(x,\cbar)\in q$ such that
$p\cup\{\theta(x,\cbar)\}$ forks over $M$.
As $p\perp\ID$, $p^{(n)}\perp\ID$, so the $\ID_\perp$-minimality
of $p$ implies $\tp(a/M\cbar)$ is $\ID$-small.  

But, since $p$ is foreign
to $\ID$,   $\tp(b/M\cbar)$, which is a nonforking extension of $p$
would be orthogonal to $q$ by Lemma~\ref{charforeign}(2).  In particular,
$\tp(b/M\cbar a)$ would not fork over $M\cbar$.
\endproof

The following easy `transfer result' will be used in the 
subsequent sections.

\begin{Lemma}  \label{transfer}
Assume that $B$ is algebraically closed, $p=\tp(a/B)$ is foreign
to $\ID$, $q=\tp(b/B)$, and $b\in\acl(Ba)\setminus B$.
Then $q$ is foreign to $\ID$.  If, in addition, $p$ is $\ID$-minimal
($\ID_\perp$-minimal) then $q$ is $\ID$-minimal ($\ID_\perp$-minimal)
as well.
\end{Lemma}

\bp  If $q$ were not foreign to $\ID$, then by Lemma~\ref{charforeign}(4)
there is $c\in\dcl(Bb)\setminus B$ such that $\tp(c/B)$ is $\ID$-small.
Since $\tp(c/B)$ is not algebraic it is not orthogonal to $p$,
which, via Lemma~\ref{charforeign}(2), contradicts $p$ being foreign to $\ID$.
Thus $q\perp\ID$.
   
Next, suppose that $p$ is $\ID$-minimal.  Since $p\not\perp q$ and $p\perp \ID$,
$q$ cannot be $\ID$-small.  To see that $q$ is $\ID$-minimal, choose $C\supseteq B$
such that $\tp(b/C)$ forks over $B$.  Then $\tp(a/C)$ forks over $B$, so $\tp(a/C)$
is $\ID$-small.  Thus $\tp(b/C)$ is $\ID$-small by Lemma~\ref{iterate}.
\endproof

\section{Chains and witnessing groups}

Throughout this section $\ID$ always denotes a regular ideal.

\begin{Definition} {\em  We say {\em $A$ is an $\ID$-subset of $B$,\/}
written $A\subseteq_{\ID} B$, if $A\subseteq B$ and $\stp(b/A)\perp\ID$
for every finite tuple $b$ from $B$.
When $M$ and $N$ are models we 
write $M\preceq_{\ID} N$ when both $M\preceq N$ and
$M\subseteq_{\ID} N$.
A set $A$ is {\em $\ID$-full\/} if $A\subseteq_{\ID} M$ for some
(equivalently for every) a-prime model $M$ over $A$.
}
\end{Definition}

\begin{Lemma} \label{ext} Let $\ID$ be any regular ideal and assume $M$ is a-saturated.
\begin{enumerate}
\item  If $M\preceq N$ are models then $M\preceq_{\ID} N$ if and only if
$\phi(N)=\phi(M)$ for all $\phi\in L(M)\cap\ID$.
\item  If $M\subseteq_{\ID} A$, then $M\preceq_{\ID} M[A]$, where
$M[A]$ is any a-prime model over $M\cup A$.
\end{enumerate}
\end{Lemma}

\bp  (1)  First suppose $M\preceq_{\ID} N$ and choose $\phi\in L(M)\cap\ID$.
If $c\in\phi(N)$ then $\tp(c/N)$ is $\ID$-small.  If $\tp(c/M)$ were not algebraic,
it would be nonorthogonal to an $\ID$-small type, contradicting $\tp(c/M)\perp\ID$.
So $\tp(c/M)$ is algebraic, hence $c\in\phi(M)$.  Conversely, if there
were $c\in N$ such that $\tp(c/M)\not\perp \ID$, then by Lemma~\ref{charforeign}(4)
there is $c'\in\dcl(Mc)\setminus M$ such that $\tp(c'/M)$ is $\ID$-small.
Then $\phi(N)\neq \phi(M)$ for any $\phi\in\tp(c'/M)\cap\ID$.

(2)  Recall that because $M$ is a-saturated,
$M[A]$ is dominated by $A$ over $M$.
Choose any tuple $c$ from $M[A]$.  If $\tp(c/M)$ were not foreign to $\ID$,
then as $M$ is a-saturated, there is an $\ID$-small type $q\in S(M)$ such that
$\tp(c/M)\not\perp q$, hence $\tp(c/M)$ is not almost orthogonal to $q$.
Since $c$ is dominated by $A$ over $M$, there is $a$ from $A$ such that
$\tp(a/M)$ is not almost orthogonal to $q$, which contradicts $M\subseteq_{\ID} A$.
\endproof

\begin{Definition} {\em 
A {\em saturated chain\/} is an elementary chain
$\<M_\alpha:\alpha<\delta\>$
of a-saturated models in which $M_{\alpha+1}$ realizes every complete
type over $M_\alpha$ for each $\alpha<\delta$.
An {\em $\ID$-chain\/} is a sequence 
$\<M_\alpha:\alpha<\delta\>$ of a-saturated models
such that $M_\alpha\preceq_{\ID} M_\beta$ for all 
$\alpha<\beta<\delta$ and $M_{\alpha+1}$
realizes every type over $M_\alpha$ foreign to $\ID$.
A chain (of either kind) is {\em $\ID$-full\/} if 
the union
$\bigcup_{\alpha<\delta} M_\alpha$ is an $\ID$-full set.
}
\end{Definition}

In general, a saturated chain need not be $\ID$-full.  However,
if $\ID$ is either
the  ideal of algebraic formulas or  superstable formulas (both of which are
regular),
then any a-saturated chain is $\ID$-full, since types are based on
finite sets.  A more complete explanation of this is given in the proof of
Lemma~\ref{ri}.
By contrast, the
following Lemma demonstrates that $\ID$-chains are always $\ID$-full.

\begin{Lemma} \label{fullness}
Every $\ID$-chain is full.  That is,
if $\<M_\alpha:\alpha<\delta\>$ is an $\ID$-chain, $\delta$ is
a nonzero limit ordinal, and $M_\delta$ is
a-prime over $\bigcup_{\alpha<\delta} M_\alpha$, then $M_\alpha\preceq_{\ID} M_\delta$
for all $\alpha<\delta$.  
\end{Lemma}

\bp  By the characterization of $M\preceq_{\ID} N$
given by Lemma~\ref{ext}(1), the first sentence follows from the second.
So fix an $\ID$-chain $\<M_\alpha:\alpha<\delta\>$.
Let $N=\bigcup_{\alpha<\delta} M_\alpha$ and let $M_\delta$ be a-prime over $N$.
Fix any $\alpha<\delta$.  Since $M_\alpha\subseteq_{\ID} M_\beta$ for all 
$\alpha<\beta<\delta$, $M_\alpha\subseteq_{\ID} N$, so $M_\alpha\preceq_{\ID} M_\delta$
by Lemma~\ref{ext}(2).
\endproof

\begin{Definition} {\em  A formula $\theta$ is {\em weakly $\ID$-minimal
(weakly $\ID_\perp$-minimal)\/}
if $\{\theta\}$ is $\ID$-minimal ($\ID_\perp$-minimal).
}
\end{Definition}

We now state offer two complementary propositions.  
The main point of both is that they produce regular types
that are `close' to a given regular ideal.
The advantage of (1) is that one obtains  $\ID$-minimality
at the cost of requiring the chain to be $\ID$-full.  In
(2) the fullness condition is automatically satisfied by Lemma~\ref{fullness},
but one only gets $\ID_\perp$-minimality.

\begin{Proposition}  \label{split}  Fix a regular ideal $\ID$, a countable,
stable theory $T$, and an $\ID$-large formula $\phi$.
\begin{enumerate}
\item
{\bf Either} there is a weakly $\ID$-minimal formula 
$\psi\vdash\phi$
{\bf or} for every $\ID$-full saturated chain $\<M_n:n\in\omega\>$ 
with $\phi\in L(M_0)$, there is an $\aleph_1$-isolated,
$\ID$-minimal $p\in S(\bigcup_n M_n)$
with $\phi\in p$ and $p\perp\ID$.  
\item
{\bf Either} there is a weakly $\ID_\perp$-minimal formula 
$\psi\vdash\phi$
{\bf or} for every $\ID$-chain $\<M_n:n\in\omega\>$ 
with $\phi\in L(M_0)$, there is an $\aleph_1$-isolated,
$\ID_\perp$-minimal $p\in S(\bigcup_n M_n)$
with $\phi\in p$ and $p\perp\ID$.  
\end{enumerate}
Moreover, in either of the two `second cases' the type $p$ is regular.
\end{Proposition}

\bp  Assume that there is no weakly $\ID$-minimal $\psi\vdash\phi$.
Fix 
an $\ID$-full  saturated chain $\<M_n:n\in\omega\>$
with $\phi\in L(M_0)$, let $N=\bigcup_{n\in\omega} M_n$,
and let $M_\omega$ be $\aleph_1$-prime over $N$. Let
$\Delta_0\subseteq\Delta_1\subseteq\dots$ be finite sets of
formulas with $L=\bigcup_{n\in\omega} \Delta_n$.  
We inductively construct a sequence $\<\phi_n:n\in\omega\>$
of $\ID$-large formulas as follows:
Let $\phi_0$ be our given $\phi$.
Given $\phi_n\vdash\phi_0$ that is an $\ID$-large $L(M_n)$-formula
$$A_n=\{\psi\in L(M_{n+1}): \psi\vdash \phi_n,
\ \psi\ \hbox{is}\ \ID\hbox{-large and
forks over}\ M_n\}.$$
As $M_{n+1}$ realizes every type over $M_n$ foreign to $\ID$
and $\phi_n$ is not
weakly $\ID$-minimal, $A_n$ is nonempty.
Choose $\phi_{n+1}\in A_n$ so as to minimize
$R(\psi,\Delta_n,2)$.
Let $\Gamma=\{\phi_n:n\in\omega\}$.  We first argue that
$\Gamma$ has a unique extension to a complete type in $S(N)$.

\medskip

{\bf Claim.}  $\Gamma\vdash\neg\psi(x,b)$ for all 
$\psi(x,b)\in \ID\cap L(N)$.

\medskip

\bp    
If the Claim were to fail, then $\Gamma\cup\{\psi(x,b)\}$ 
would be consistent,
hence would be realized in $M_\omega$, say by an element $c$.
As the chain is $\ID$-full, $c\in N$.
Choose $n$ such that
$b,c\in M_n$.  But $\phi_{n+1}$ was chosen to fork over
$M_n$, yet is realized in $M_n$, which is impossible.
\endproof

Now let $\psi(x,b)$ be any $L(N)$-formula.
Choose $n$ such that $\psi(x,y)\in\Delta_n$.  
As $\phi_{n+1}$ was chosen to be of minimal R(--,$\Delta_n,2)$-rank,
it is not possible for both $\phi_{n+1}\wedge\psi(x,b)$ and
$\phi_{n+1}\wedge\neg\psi(x,b)$ to be in $A_n$.  As
$\ID$ is an ideal, at least one of the two of them is $\ID$-large,
so is an element of $A_n$, thus the other one is $\ID$-small or
inconsistent.  Using the Claim, either $\Gamma\vdash\psi(x,b)$
or $\Gamma\vdash\neg\psi(x,b)$.
Thus $\Gamma$ implies a complete type in $S(N)$, which we
call $p$.  

By construction $p$ is $\aleph_1$-isolated 
and is $\ID$-large by the Claim.  Since $M_\omega$ is $\aleph_1$-saturated and $p$ is
$\aleph_1$-isolated, there is a realization $c$ of $p$ in $M_\omega$.
If $p$ were not foreign to $\ID$ then by Lemma~\ref{charforeign}(4) there would be
$c'\in\dcl(Nc)\setminus N$ with $c'/N$ $\ID$-small, directly
contradicting $\ID$-fullness.

It remains to show that any forking extension of $p$ is $\ID$-small.
let $\theta(x,a^*)$ be any $L(\C)$-formula such that $p\cup\theta(x,a^*)$
forks over $M_\omega$.  Then for some $n$,
$\phi_{n+1}\wedge\theta(x,a^*)$ $\Delta_n$-forks over $M_n$.
As $M_{n+1}$ realizes all types over $M_n$
there is $a'\in M_{n+1}$ such that $\tp(a'/M_n)=\tp(a^*/M_n)$.
But then $\phi_{n+1}\wedge\theta(x,a')$ $\Delta_n$-forks over
$M_n$, contradicting the 
minimality of R(--,$\Delta_n,2)$ rank of $\phi_{n+1}$.

As for (2) assume that there is no $\ID_\perp$-minimal formula
implying $\phi$.  Choose an $\ID$-chain $\<M_n:n\in\omega\>$,
which is automatically $\ID$-full by Lemma~\ref{fullness}.
The definition of $\{A_n\}_{n\in\omega}$ and the constructions of
$\Gamma$ and $p$ remain the same.  All that is affected is that
in the final paragraph, as we only need to establish $\ID_\perp$-minimality,
one chooses a formula $\theta(x,a^*)$ with $\tp(a^*/N)\perp\ID$.
By Lemma~\ref{charforeign}(4) this implies $\tp(a^*/M_n)\perp\ID$
for all $n\in\omega$, so choosing $n$ as above, one obtains
$a'\in M_{n+1}$ satisfying $\tp(a'/M_n)=\tp(a^*/M_n)$ and
a similar contradiction is obtained.

In both cases, the regularity of $p$ follows immediately from Lemma~\ref{regular}.
\endproof

Recall that a stable theory has {\em NDIDIP\/} if for every elementary
chain $\<M_n:n\in\omega\>$ of models, every type that is nonorthogonal
to some a-prime model over $\bigcup_{n\in\omega} M_n$ is nonorthogonal to some
$M_n$.  Relationships between NDIDIP and NDOP are explored in \cite{LSh}.

\begin{Proposition}  \label{group}
Fix a countable, stable theory $T$ with NDIDIP and a regular ideal $\ID$
such that the formula `$x=x$'$\not\in\ID$.
\begin{enumerate}
\item   If there is an
an $\ID$-full,  saturated chain
$\<M_n:n\in\omega\>$, but there is no weakly $\ID$-minimal formula then
there is an abelian group witness to unsuperstability, where in addition
the generic type of the intersection
is both $\ID$-minimal and foreign to $\ID$.
\item If there is no weakly $\ID_\perp$-minimal formula then
there is an abelian group witness to unsuperstability where the generic
type of the intersection is $\ID_\perp$-minimal and foreign to $\ID$.
\end{enumerate}
\end{Proposition}

\bp  (1) Fix an $\ID$-full, saturated chain $\<M_n:n\in\omega\>$ and let 
$N=\bigcup_{n\in\omega} M_n$.  
Using Proposition~\ref{split}(1)
choose $p\in S(N)$ to be $\aleph_1$-isolated, foreign to $\ID$, and $\ID$-minimal, hence
regular.  Since $T$ has NDIDIP, $p\not\perp M_n$.  Since $p$ is regular
and $M_n$ is a-saturated,  by Claim~X~1.4 of \cite{Shc}
there is a regular type $r_0\in S(M_n)$ nonorthogonal to $p$.  
Let $r$ denote the nonforking
extension of $r_0$ to $N$.  
As $p$ and $r$  are nonorthogonal  there is an integer
$m$ such that $p^{(m)}$ is not almost
orthogonal to $r^{(\omega)}$.  
Since $p$ is $\aleph_1$-isolated and $M_n$ is
a-saturated, $Na$ is dominated by $N$ over $M_n$ for any $a$ realizing $p$.
Thus $p^{(1)}$  is not almost orthogonal to $r^{(\omega)}$ over $N$. 
Choose $k\ge 1$ maximal such that 
$p^{(k)}$ is almost orthogonal to $r^{(\omega)}$ over
$N$ and choose $\cbar$ realizing $p^{(k)}$. Let 
$B=\acl(N\cbar)$ and choose a realization
$a$ of the nonforking extension of $p$ to $B$.   
 
By Theorem~1 of \cite{Hr2}, there is $b\in\dcl(Ba)\setminus B$
and a type-definable, connected group 
$A$ with a regular generic type $q$ 
(so $A$ is abelian by Poizat's theorem~\cite{Poizat})
and a definable regular, transitive action of $A$ on $p_1(\C)$, where
$p_1=\tp(b/B)$.  By Lemma~\ref{transfer} the type $p_1$ and hence $q$ are both
foreign to $\ID$ and $\ID$-minimal.  By Theorem~2 of \cite{Hr3} there is
a definable supergroup $A_0\supseteq A$.  By an easy compactness argument we may
assume $A_0$ is abelian as well.  Furthermore, by iterating Theorem~2 of
\cite{Hr3} we obtain a descending sequence $\<A_n:n\in\omega\>$ of subgroups of $A_0$
with $A=\bigcap_{n\in\omega} A_n$.  

Thus far we have not guaranteed that $A_{n+1}$ has infinite index in $A_n$.  
In order to show that there is a subsequence of the $A_n$'s with this property
and thereby complete the proof of the Proposition,
it suffices to prove the following claim:

\medskip

{\bf Claim}  For every $n\in\omega$ there is $m\ge n$ such that $[A_n:A_m]$ is infinite.

\medskip

\bp  By symmetry it suffices to show this for $n=0$.  Assume that this were not the
case, i.e., that $[A_0,A_m]$ is finite for each $m$.  Then $A$ has bounded index in $A_0$.
We will obtain a contradiction by showing that the definable set $A_0$ is
weakly $\ID$-minimal.
First, since $q$ is $\ID$-large, the formula defining $A_0$
is $\ID$-large as well.
Let $\phi(x,e)$ be any forking extension of the formula defining $A_0$
and let $E\subseteq A_0$ be the set of realizations of $\phi(x,e)$.
Let $\{C_i:i<2^\kappa\le 2^{\aleph_0}\}$ 
enumerate the $A$-cosets of $A_0$.
For each $i$, $E\cap C_i$ is a forking extension of $C_i$.
Since every $C_i$ is a translate of $A$ whose generic type is
$\ID$-minimal, this implies that $E\cap C_i$ is $\ID$-small for
each $i$.  Thus $\phi(x,e)\in\ID$ by compactness (and the fact
that $\ID$ is an ideal).  Thus, the formula defining $A_0$ 
is weakly $\ID$-minimal, contradiction.

The proof of (2) is identical, choosing an $\ID$-chain satisfying
the hypotheses and using Proposition~\ref{split}(2)
in place of \ref{split}(1).
\endproof

\section{Applications}  \label{app}

Our first application gives a `trichotomy' for strictly stable
theories in a countable language.  It uses the
ideal of superstable formulas.  Let $\Rinfty$ denote the ideal of

\begin{Definition}  \label{ss}  {\em 
$\Rinfty$ denotes the ideal of 
superstable formulas (i.e., all formulas $\phi$  with $R^\infty(\phi)<\infty$).
}
\end{Definition}

Equivalently, $\phi\in\Rinfty$ if and only if for all  cardinals
$\kappa\ge 2^{|T|}$, for any model $M$ of size $\kappa$
containing the parameters of $\phi$, there are at most
$\kappa$ complete types over $M$ extending $\phi$.

\begin{Lemma}  \label{ri}
$\Rinfty$ is a regular ideal,
any elementary chain $\<M_n:n\in\omega\>$ of a-saturated models is 
$\Rinfty$-full, and there are no weakly $\Rinfty$-minimal formulas.
\end{Lemma}

\bp  Invariance under automorphisms of $\C$ is clear and $\Rinfty$ being
an ideal follows by counting types.  To show regularity,
choose $\psi(y)\in\Rinfty$ and $\theta(x,y)\in L(\C)$ such that 
$\theta(x,b)\in\Rinfty$ for every $b$ realizing $\psi$.
Choose $\kappa\ge 2^{|T|}$  and a model $M$ of size $\kappa$
containing the hidden parameters of both $\psi$ and $\theta$.
Then there are at most $\kappa$ types $p(x,y)\in S(M)$ extending
$\theta(x,y)\wedge\psi(y)$, so the projection 
$\exists y(\theta(x,y)\wedge\psi(y))\in\Rinfty$ as only $\kappa$
types $q(x)\in S(M)$ extend it.

To establish fullness, fix an elementary chain $\<M_n:n\in\omega\>$
of a-saturated models.
Let $N=\bigcup_{n\in\omega} M_n$ and choose an a-prime model $M_\omega$
over $N$. Because of Lemma~\ref{charforeign}(4), in order to show
that $N\subseteq_{\ID} M_\omega$ it suffices to show that
no element of  
$c\in M_\omega\setminus N$ is $\Rinfty$-small.  So choose
$c\in M_\omega$ such that $\tp(c/N)$ is $\ID$-small and we will show
that $c\in N$.
On one hand, since $\tp(c/N)$ contains a superstable formula
there is a finite $n$ such that $\tp(c/N)$ is based on $M_n$.
On the other hand, since $M_\omega$ is a-prime over $N$,
$\tp(c/N)$ is a-isolated.
Thus $\tp(c/M_n)$ is a-isolated as well (see
e.g., Theorem~IV~4.3(1) of \cite{Shc}).
Since $M_n$ is a-saturated, this implies $c\in M_n\subseteq N$.

To show that there are no weakly $\Rinfty$-minimal formulas,
suppose that a formula $\phi$ has the property that any forking
extension of $\phi$ is $\Rinfty$-small.
We will show that $\phi\in\Rinfty$ by counting types.
Fix a cardinal $\kappa\ge 2^{|T|}$ and  a model $M$ of size $\kappa$
containing the parameters of $\phi$.  Let $M_0\preceq M$
have size $|T|$ that also contains the parameters containing $\phi$.
It suffices to show that every $p\in S(M_0)$ extending $\phi$
has at most $\kappa$ extensions to types in $S(M)$.  
Clearly, there is a unique nonforking extension of $p$
and any forking extension of $p$ contains an $L(M)$-formula
witnessing the forking.  Each such forking formula $\psi\in\Rinfty$,
so there are at most $\kappa$ $q\in S(M)$ extending $\psi$.
So, since the total number of $\psi\in L(M)$ is at most $\kappa$,
$p$ has at most $\kappa$ extensions to types in $S(M)$.
\endproof

\begin{Theorem}  
Let $T$ be a stable, unsuperstable theory in a countable language.
Then at least one of the following three conditions occurs:
\begin{enumerate}
\item $T$ has the dimensional order property (DOP); or
\item $T$ has NDOP, but is deep (i.e., there is a sequence 
$\<M_n:n\in\omega\>$ such that $\tp(M_{n+1}/M_n)\perp M_{n-1}$
for all $n\ge 1$); or
\item There is an abelian group witness to unsuperstability
(see Definition~\ref{witness})
in which the generic type of the intersection is both
$\Rinfty$-minimal and foreign to $\Rinfty$.
\end{enumerate}
\end{Theorem}

\bp  
To begin, Corollary 1.12 of \cite{LSh} asserts
that  any such  theory $T$ has NDIDIP.
Since $T$ is not superstable the formula `$x=x$'$\not\in\Rinfty$.
As well, by Lemma~\ref{ri}  there are no weakly
$\Rinfty$-minimal formulas, so
Proposition~\ref{group}(1) asserts that an abelian group
witness to unsuperstability exists, whose generic type is
regular and 
both $\Rinfty$-minimal and foreign to $\Rinfty$.
\endproof

Our second application comes from an attempt to solve the
`Main Gap for $\aleph_1$-saturated models.'  As in the previous
theorem, the relevant setting is where a countable theory $T$
is stable, unsuperstable, with NDOP, and is shallow.
The main open question is whether, for such a theory
every nonalgebraic type $r$ is nonorthogonal to a regular type.
The following result sheds some light on this issue.
In order to analyze this problem, fix a nonalgebraic, stationary type $r$ over the
empty set.  Let 
$$\ID_r=\{\phi\in L(\C):r\perp \phi\}$$
Verifying that $\ID_r$ is an invariant ideal is straightforward.
To see that it is a regular ideal, fix $L(\C)$-formulas $\psi(y)\in\ID_r$
and $\theta(x,y)$ such that $\theta(x,b)\in\ID_r$ for every $b$ realizing
$\psi$.  Choose an a-saturated model $M$ containing the parameters of 
$\psi$ and $\theta$, pick a realization $c$ of the nonforking
extension of $r$ to $M$, and let $M[c]$ be 
any a-prime model over $Mc$.
To show that $\phi(x):=\exists y(\theta(x,y)\wedge\psi(y))\perp r$ it
suffices to prove that any realization of $\phi$ in $M[c]$ is contained
in $M$.  So choose any $a\in\phi(M[c])$.  Choose $b\in M[c]$
such that $\theta(a,b)\wedge \psi(b)$ holds.  Since $r\perp\psi$, $b\in M$.
But then $\theta(x,b)$ is over $M$ and is $\perp r$, so $a\in M$ as well.
Thus $\ID_r$ is a regular ideal.

\begin{Theorem}  \label{rperp}  Assume that a countable 
theory $T$ is stable, unsuperstable,
has NDOP, and is shallow.  If a nonalgebraic, 
stationary type $r$ is orthogonal to every regular
type, then there is an abelian group witness to unsuperstability in which the
generic type of the intersection $A=\bigcup_n A_n$ is both $(\ID_r)_\perp$-minimal
and foreign to $\ID_r$.
\end{Theorem}

\bp  Fix such a type $r$.  By naming constants we may assume that $r$ is over
the empty set.  Note that any formula $\phi\in r$ is not an element of $\ID_r$,
so `$x=x$'$\not\in\ID_r$.

\medskip

{\bf Claim.}  There is no weakly $(\ID_r)_\perp$-minimal formula.

\medskip 

\bp  Assume that $\phi$ were $(\ID_r)_\perp$-minimal.  We  construct
a regular type $p\not\perp r$ as follows:  Choose an a-saturated model
$M$ containing the parameters in $\phi$, pick a realization $c$ of
the nonforking extension of $r$ to $M$, and choose
an a-prime model $M[c]$ over $Mc$.  Since $\phi$ is $\ID_r$-large
we can find an  $a\in M[c]\setminus M$ realizing $\phi$.  Choose such an $a$
and let $p=\tp(a/M)$.  Clearly, $p\not\perp r$.  To see that $p$ is regular,
first note that $p$ is $(\ID_r)$-minimal
since $p$ is $\ID_r$-large and extends $\phi$.
As well, $p$ is foreign to $\ID_r$, since if it were not, then by
Lemma~\ref{charforeign}(4) there would be $b\in\dcl(Ma)$ with $\tp(b/M)$
$\ID_r$-small.  But then $tp(c/Mb)$ would fork over $M$, implying that
$r$ is nonorthogonal to an $\ID_r$-small type, which is a contradiction.
So $p$ is $(\ID_r)$-minimal and foreign to $\ID_r$, hence is regular
by Lemma~\ref{regular}.
\endproof

The theorem now follows immediately from Proposition~\ref{split}(2).
\endproof

\end{document}